\documentclass[12pt]{article}

\usepackage{epsfig}
\usepackage{amssymb}
\usepackage{amsmath}
\usepackage{subfigure}

\def\shadowbox{\hbox{\rule[-0.0ex]{0.1ex}{1.2ex}%
\hspace{-0.1ex}\rule[-0.0ex]{1.2ex}{0.1ex}%
\hspace{0.0ex}\rule[-0.0ex]{0.1ex}{1.2ex}\hspace{-1.3ex}%
\rule[1.15ex]{1.25ex}{0.1ex}\hspace{-0.0ex}\rule[-0.25ex]{0.3ex}{1.1ex}%
\hspace{-1.2ex}\rule[-0.25ex]{1.1ex}{0.25ex}}}
\def\qed{\ifmmode \hbox{\hfill\shadowbox}
     \else \hphantom{x}\hfill\shadowbox \fi}

\newtheorem{theorem}{Theorem}[section]
\newtheorem{lemma}[theorem]{Lemma}

\def\Cst{{\mathbb C}}
\def\Nst{{\mathbb N}}
\def\Rst{{\mathbb R}}
\def\Rdst{{\Rst^d}}
\def\Tst{{\mathbb T}}
\def\Zst{{\mathbb Z}}
\def\Lsp{{\boldsymbol L}}
\def\Ltsp{{\Lsp^2}}
\def\LtR{{\Ltsp(\Rst)}}

\def\lsp{{\boldsymbol\ell}}

\def\ltZ{{\lsp^2(\Zst})}
\def\Hsp{{\boldsymbol H}}
\def\sinc{\mbox{sinc\/}}
\def\sinco{\mbox{sinc\/}_{\Omega}}
\def\diag{\mbox{diag\/}}
\def\eps{\varepsilon}
\def\ord{{\cal O }}
\def\supp{\mbox{supp\/}}
\def\cond{\mbox{cond\/}}

\def\circ{\operatorname{circ}}
\def\range{\operatorname{\mathcal R}}

\def\thresh{{\tau}}
\def\Bsp{{\boldsymbol B}}
\def\BO{{\Bsp}}
\def\BOM{{\Bsp_{\Omega}}}
\def\Psp{{\boldsymbol P}}
\def\PM{{\Psp_{\! M}}}
\def\PMd{{\Psp^{d}_{\! M}}}

\def\sinco{{\sinc}}
\def\fk{{f_j}}
\def\fdk{{\gamma_j}}
\def\Pn{{P_{n}}}
\def\Ri{{R^{-1}}}
\def\Rp{{R^{+}}}
\def\Rn{{R_{n}}}
\def\Rnp{{R_{n}^{+}}}
\def\Rnpe{{R_{n}^{+,\thresh}}}
\def\Rni{{R_{n}^{-1}}}

\def\cn{{c^{(n)}}}
\def\cnd{{c^{(n,\delta)}}}
\def\bn{{b^{(n)}}}
\def\bnd{{b^{(n,\delta)}}}
\def\fn{{f^{(n)}}}

\def\sincframe{{\{\sinco(\cdot - t_j)\}_{j \in \Zst}}}
\def\ntoinf{{n \rightarrow \infty }}
\def\toinf{{\rightarrow \infty }}
\def\plsp{{p_{M}}}
\def\alsp{{a_{M}}}
\def\alspk{{(\alsp)_{k}}}

\begin{document}

\title{\bf Numerical Analysis of the Non-uniform Sampling Problem}
\author{Thomas Strohmer\thanks{Department of Mathematics, University
of California, Davis, CA-95616; strohmer@math.ucdavis.edu.
\hspace*{5.4mm} The author was supported by NSF DMS grant 9973373.}}

\date{}
\maketitle

\begin{abstract}

We give an overview of recent developments in the problem of reconstructing 
a band-limited signal from non-uniform sampling from a numerical analysis 
view point. It is shown that the appropriate design of the finite-dimensional 
model plays a key role in the numerical solution of the non-uniform
sampling problem.
In the one approach (often proposed in the literature) the finite-dimensional 
model leads to an ill-posed problem even in very simple situations. The other
approach that we consider leads to a well-posed problem that preserves 
important structural properties of the original infinite-dimensional problem
and gives rise to efficient numerical algorithms. Furthermore a fast
multilevel algorithm is presented that can reconstruct
signals of unknown bandwidth from noisy non-uniformly spaced samples.
We also discuss the design of efficient regularization methods for
ill-conditioned reconstruction problems.
Numerical examples from spectroscopy and exploration geophysics
demonstrate the performance of the proposed methods. 
\end{abstract}
\mbox{Key words:} non-uniform sampling, band-limited functions, frames,
regularization, signal reconstruction, multi-level method.

\section{Introduction}
\label{s:intro}

The problem of reconstructing a signal $f$ from non-uniformly spaced
measurements $f(t_j)$ arises in areas as diverse as geophysics, medical 
imaging, communication engineering, and astronomy. 
A successful reconstruction of $f$ from its samples $f(t_j)$ requires
a priori information about the signal, otherwise the reconstruction
problem is ill-posed. This a priori information can often be obtained from 
physical properties of the process generating the signal.
In many of the aforementioned applications the signal can be assumed to 
be (essentially) band-limited.

Recall that a signal (function) is band-limited with bandwidth $\Omega$
if it belongs to the space $\BOM$, given by
\begin{equation}
\BOM = \left\{ f \in \LtR : \hat{f}(\omega) = 0
          \,\,\text{for}\,\,|\omega|>\Omega \right\}\,,
\label{bandlim}
\end{equation}
where $\hat{f}$ is the Fourier transform of $f$ defined by 
$$\hat{f}(\omega) 
= \int \limits_{-\infty}^{+\infty} f(t) e^{-2 \pi i \omega t} \, dt\,.$$
For convenience and without loss of generality we restrict our attention
to the case $\Omega = \frac{1}{2}$, since any other bandwidth can be
reduced to this case by a simple dilation. Therefore we will
henceforth use the symbol $\BO$ for the space of band-limited signals.
 
It is now more than 50 years ago that Shannon published his celebrated 
sampling theorem~\cite{Sha48}. His theorem implies that any signal 
$f \in \BO$ can be reconstructed from
its regularly spaced samples $\{f(n)\}_{n \in \Zst}$ by
\begin{equation}
f(t) = \sum_{n \in \Zst} f(n)
    \frac{\sin \pi(t - n)}{\pi(t-n)}\,.
\label{shannon}
\end{equation}

In practice however we seldom enjoy the luxury of equally spaced samples.
The solution of the nonuniform sampling problem poses much more difficulties,
the crucial questions being: 
\begin{itemize}
\vspace*{-1mm}
\setlength{\itemsep}{-0.5ex}
\setlength{\parsep}{-0.5ex}
\item Under which conditions is a signal $f \in \BO$ uniquely defined 
      by its samples $\{f(t_j)\}_{j \in \Zst}$?
\item How can $f$ be stably reconstructed from its samples $f(t_j)$?
\end{itemize}

These questions have led to a vast literature on nonuniform sampling theory
with deep mathematical contributions see~\cite{DS52,Lan67,BM67,BSS88,FG94}
to mention only a few. There is also no lack of methods claiming 
to efficiently reconstruct a function from its 
samples~\cite{Yen56,YT67,Ben92,FGS95,Win92,Mar93a,FG94}.
These numerical methods naturally have to operate in a finite-dimensional
model, whereas theoretical results are usually derived for the
infinite-dimensional space $\BO$. From a numerical point of view the
``reconstruction'' of a bandlimited signal $f$ from a finite number of samples
$\{f(t_j)\}_{j=1}^{r}$ amounts to computing an approximation to $f$ 
(or $\hat{f}$) at sufficiently dense (regularly) spaced grid points in an 
interval $(t_1, t_r)$.

Hence in order to obtain a  ``complete'' solution of the sampling problem
following questions have to be answered: 
\begin{itemize}
\vspace*{-1mm}
\setlength{\itemsep}{-0.5ex}
\setlength{\parsep}{-0.5ex}
\item Does the approximation computed within the finite-dimensional model 
actually converge to the original signal $f$, when the dimension of
the model approaches infinity?
\item Does the finite-dimensional model give rise to fast and stable
numerical algorithms?
\end{itemize}

These are the questions that we have in mind, when presenting an overview
on recent advances and new results in the nonuniform sampling problem
from a numerical analysis view point.

In Section~\ref{ss:truncated} it is demonstrated that the celebrated frame
approach does only lead to fast and stable numerical methods when the
finite-dimensional model is carefully designed. The approach usually
proposed in the literature leads to an ill-posed problem even in 
very simple situations. We discuss several methods to stabilize the 
reconstruction algorithm in this case.
In Section~\ref{ss:trigpol} we derive an alternative finite-dimensional model,
based on trigonometric polynomials. This approach leads to a well-posed problem
that preserves important structural properties of the original 
infinite-dimensional problem and gives rise to efficient numerical algorithms. 
Section~\ref{s:numeric} describes how this approach can be modified
in order to reconstruct band-limited signals for the in practice
very important case when the bandwidth of the signal is not known.
Furthermore we present regularization techniques for ill-conditioned 
sampling problems. Finally Section~\ref{s:applications} 
contains numerical experiments from spectroscopy and geophysics.

Before we proceed we introduce some notation that will be used throughout
the paper.
If not otherwise mentioned $\|h\|$ always denotes the $\LtR$-norm
($\ltZ$-norm) of a function (vector).
For operators (matrices) $\|T\|$ is the standard operator (matrix) norm.
The condition number of an invertible operator $T$ is defined by 
$\kappa (A) = \|A\| \|A^{-1}\|$ and the spectrum of $T$ is $\sigma (T)$.
$I$ denotes the identity operator.

\subsection{Nonuniform sampling, frames, and numerical algorithms} 
\label{s:theory}

The concept of frames is an excellent tool to study nonuniform sampling 
problems~\cite{Fei89,BH90,Ben92,Hig96,FG94,Zay93}. The frame approach has 
the advantage that it gives rise to deep theoretical results and also to
the construction of efficient numerical algorithms -- {\em if} (and this 
point is often ignored in the literature) the finite-dimensional model is 
properly designed.

Following Duffin and Schaeffer~\cite{DS52}, a family $\{\fk\}_{j \in \Zst}$ 
in a separable Hilbert space $\Hsp$ is said to be a frame for $\Hsp$, if
there exist constants (the {\em frame bounds}) $A,B>0$ such that
\begin{equation}
\label{framedef}
A \|f\|^2 \le \sum_{j} |\langle f, \fk \rangle|^2 \le B \|f\|^2 \,,
\qquad \forall f \in \Hsp.
\end{equation}
We define the {\em analysis operator} $T$ by
\begin{equation}
T: f \in \Hsp \rightarrow Ff = \{ \langle f, \fk \rangle\}_{j \in \Zst}\,,
\label{frameanal}
\end{equation}
and the {\em synthesis operator}, which is just the adjoint operator of
$T$, by
\begin{equation}
T^{\ast}: c \in \ltZ \rightarrow
T^{\ast} c = \sum_{j} c_{j} \fk\,.
\label{framesyn}
\end{equation}
The {\em frame operator} $S$ is defined by $S = T^{\ast} T$, hence
$Sf = \sum_{j} \langle f, \fk \rangle \fk$. $S$ is bounded
by $A I \le S \le B I$ and hence invertible on $\Hsp$.

We will also make use of the operator $T T^{\ast}$ in form of
its Gram matrix representation $R: \ltZ \rightarrow \ltZ $ with entries
$R_{j,l} = \langle f_j, f_l \rangle$. On $\range (T) = \range (R)$
the matrix $R$ is bounded by $A I \le R \le B I$ and invertible.
On $\ltZ$ this inverse extends to the {\em Moore-Penrose inverse} or
pseudo-inverse $R^{+}$ (cf.~\cite{EHN96}).

Given a frame $\{\fk\}_{j \in \Zst}$ for $\Hsp$, any $f \in \Hsp$ can be 
expressed as
\begin{equation}
\label{frameexp}
f = \sum_{j \in \Zst} \langle f, \fk \rangle \fdk 
 = \sum_{j \in \Zst} \langle f, \fdk \rangle \fk \,,
\end{equation}
where the elements $\fdk :=S^{-1} \fk$ form the so-called dual frame
and the frame operator induced by $\fdk$ coincides with $S^{-1}$.
Hence if a set $\{\fk\}_{j \in \Zst}$ establishes a frame for $\Hsp$,
we can reconstruct any function $f \in \Hsp$ from its moments
$\langle f, \fk \rangle$.

One possibility to connect sampling theory to frame theory is by means of the
{\em sinc}-function
\begin{equation}
\sinco(t) = \frac{\sin \pi t}{\pi t}\,.
\label{sinc}
\end{equation}
Its translates give rise to a {\em reproducing kernel} for $\BO$ via
\begin{equation}
f(t) = \langle f, \sinco(\cdot - t) \rangle \quad \forall t,  f \in \BO\,.
\label{sincconv}
\end{equation}
Combining~\eqref{sincconv} with formulas~\eqref{framedef} and~\eqref{frameexp}
we obtain following well-known result~\cite{Fei89,BH90}.
\begin{theorem}
If the set $\sincframe$ is a frame for $\BO$, then the function 
$f \in \BO$ is uniquely defined by the sampling set $\{f(t_j)\}_{j \in \Zst}$.
In this case we can recover $f$ from its samples by
\begin{equation}
\label{recon1}
f(t) = \sum_{j \in \Zst} f(t_j) \gamma_j \,,
\qquad \text{where}\,\,\, \gamma_j = S^{-1} \sinco(\cdot - t_j)\,,
\end{equation}
or equivalently by
\begin{equation}
\label{recon2}
f(t) = \sum_{j \in \Zst} c_j \sinco(t-t_j)\,,
\qquad \text{where}\,\,\, Rc=b\,,
\end{equation}
with $R$ being the frame Gram matrix with entries $R_{j,l}= \sinco(t_j - t_l)$
and $b=\{b_j\}=\{f(t_j)\}$. 
\end{theorem}

The challenge is now to find easy-to-verify conditions for the sampling 
points $t_j$ such that $\sincframe$ (or equivalently the exponential system 
$\{e^{2 \pi i t_j \omega}\}_{j \in \Zst}$) is a frame for $\BO$.
This is a well-traversed area (at least for one-dimensional signals),
and the reader should consult~\cite{Ben92,FG94,Hig96} for further
details and references. If not otherwise mentioned from now on we will
assume that $\sincframe$ is a frame for $\BO$.

Of course, neither of the formulas~\eqref{recon1} and~\eqref{recon2} can
be actually implemented on a computer, because both involve the solution
of an infinite-dimensional operator equation, whereas
in practice we can only compute a finite-dimensional approximation.
Although the design of a valid finite-dimensional model poses severe
mathematical challenges, this step is often neglected in theoretical
but also in numerical treatments of the nonuniform sampling problem.
We will see in the sequel that the way we design our
finite-dimensional model is crucial for the stability and efficiency
of the resulting numerical  reconstruction algorithms.

In the next two sections we describe two different approaches for obtaining 
finite-dimensional approximations to the formulas~\eqref{recon1} 
and~\eqref{recon2}. The first and more traditional approach, discussed in 
Section~\ref{ss:truncated}, applies a finite section method to 
equation~\eqref{recon2}. This approach leads to an ill-posed problem 
involving the solution of a large unstructured linear system of equations. 
The second approach, outlined in Section~\ref{ss:trigpol}, constructs a finite 
model for the operator equation in~\eqref{recon1} by means of trigonometric 
polynomials. This technique leads to a well-posed problem that is
tied to efficient numerical algorithms.

\section{Truncated frames lead to ill-posed problems} \label{ss:truncated}

According to equation~\eqref{recon2} we can reconstruct $f$ from its 
sampling values $f(t_j)$ via $f(t) = \sum_{j \in \Zst} c_j\, \sinco(t - t_j)$,
where $c=\Rp b$ with $b_j = f(t_j), j \in \Zst$. 
In order to compute a finite-dimensional approximation to 
$c = \{c_j\}_{j \in \Zst}$ we use the finite section method \cite{GF74}. 
For $x \in \ltZ$ and $n \in \Nst$ 
we define the orthogonal projection $\Pn$ by
\begin{equation}
\label{defP}
\Pn x = (\dots, 0,0, x_{-n}, x_{-n+1},\dots , x_{n-1}, x_n, 0,0, \dots)
\end{equation}
and identify the image of $\Pn$ with the space $\Cst^{2n+1}$. 
Setting $\Rn = \Pn R \Pn$ and $\bn = \Pn b$, we
obtain the $n$-th approximation $\cn$ to $c$ by solving
\begin{equation}
\Rn \cn = \bn\,.
\label{finitesec}
\end{equation}

It is clear that using the truncated frame 
$\{\sinco(\cdot - t_j)\}_{j=-n}^{n}$ in~\eqref{recon2} for an
approximate reconstruction of $f$ leads to the same system of equations.

If $\sincframe$ is an exact frame (i.e., a Riesz basis) for $\BO$ then we 
have following well-known result.
\begin{lemma}
Let $\sincframe$ be an exact frame for $\BO$ with frame bounds $A,B$
and $Rc=b$ and $\Rn \cn = \bn$ as defined above. Then $\Rni$ converges 
strongly to $\Ri$ and hence $\cn \rightarrow c$ for $\ntoinf$.
\end{lemma}
Since the proof of this result given in~\cite{Chr96b} is somewhat lengthy
we include a rather short proof here.

\begin{proof}
Note that $R$ is invertible on $\ltZ$ and $A \le R \le B$. 
Let $x \in \Cst^{2n+1}$ with $\|x\| =1$, then
$\langle \Rn x, x \rangle = \langle \Pn R \Pn x, x \rangle = 
\langle Rx,x \rangle \ge A$.
In the same way we get $\|\Rn \| \le B$,
hence the matrices $\Rn$ are invertible and uniformly bounded by 
$A \le \Rn \le B$ and 
$$\frac{1}{B} \le \Rni \le \frac{1}{A} \qquad \text{for all} \,\,n \in \Nst.$$
The Lemma of Kantorovich~\cite{RM94} yields that 
$\Rni \rightarrow \Ri$ strongly. 
\end{proof}

If $\sincframe$ is a non-exact frame for $\BO$ the situation is more 
delicate. Let us consider following situation.

\noindent
{\bf Example 1:}
Let $f \in \BO$ and let the sampling points be given by 
$t_j = \frac{j}{m}, j \in \Zst, 1 < m \in \Nst$, i.e., the signal is regularly
oversampled at $m$ times the Nyquist rate. In this case the reconstruction
of $f$ is trivial, since the set $\{\sinco(\cdot - t_j)\}_{j \in \Zst}$ is a 
tight frame with frame bounds $A=B=m$. Shannon's Sampling Theorem implies
that $f$ can be expressed as 
$f(t) = \sum_{j \in \Zst} c_j \, \sinco(t-t_j)$
where $c_j = \frac{f(t_j)}{m}$ and the numerical approximation is
obtained by truncating the summation, i.e.,
$$f_n(t) = \sum_{j =-n}^{n} \frac{f(t_j)}{m}\, \sinco(t-t_j)\,.$$

Using the truncated frame approach one finds that $R$ is a Toeplitz matrix
with entries 
$$R_{j,l}=\frac{\sin\frac{\pi}{m} (j-l)}{\frac{\pi}{m}(j-l)} 
\,, \qquad j,l \in \Zst\,,$$
in other words, $\Rn$ coincides with the prolate matrix~\cite{sle78,Var93}.
The unpleasant numerical properties of the prolate matrix are
well-documented. In particular we know that the singular values $\lambda_n$ 
of $\Rn$ cluster around $0$ and $1$ with $\log n$ singular values in the 
transition region. Since the singular values of $\Rn$ decay exponentially 
to zero the finite-dimensional reconstruction problem has become 
{\em severely ill-posed}~\cite{EHN96}, although the infinite-dimensional 
problem is ``perfectly posed'' since the frame operator satisfies $S = mI$,
where $I$ is the identity operator. 

\medskip

Of course the situation does not improve when we consider non-uniformly
spaced samples.  In this case it follows from standard linear algebra
that $\sigma(R) \subseteq \{0 \cup [A,B]\}$, or expressed in 
words, the singular values of $R$ are 
bounded away from zero. However for the truncated matrices $\Rn$ we have 
$$\sigma (\Rn)\subseteq \{(0,B]\}$$ and the smallest of the singular values
of $\Rn$ will go to zero for $\ntoinf$, see~\cite{Har98}. 

Let $A=U\Sigma V^{\ast}$ be the singular value decomposition of a matrix
$A$ with $\Sigma = \diag(\{\lambda_{k}\})$.
Then the Moore-Penrose inverse of $A$ is 
$A^+ = V \Sigma^+ U^{\ast}$, where (e.g., see~\cite{GL96})
\begin{equation}
\Sigma^{+} = \diag(\{\lambda_{k}^{+}\})\,, \quad 
\lambda_{k}^{+} = 
\begin{cases}
1/\lambda_k & \text{if}\,\, \lambda_k \neq 0, \\
0 & \text{otherwise.}
\end{cases}
\label{pinv}
\end{equation}
For $\Rn = U_n \Sigma_n V_n$ this means that the singular values close to
zero will give rise to extremely large coefficients in $\Rnp$. In
fact $\|\Rnp\| \toinf$ for $\ntoinf$ and consequently $\cn$ does not
converge to $c$. 

Practically $\|\Rnp\|$ is always bounded due
to finite precision arithmetics, but it is clear that it will lead
to meaningless results for large $n$.
If the sampling values are perturbed due to round-off error or data
error, then those error components which correspond to small
singular values $\lambda_k$ are amplified by the (then large)
factors $1/\lambda_k$. Although for a given $\Rn$ these amplifications
are theoretically bounded, they may be practically unacceptable large.

Such phenomena are well-known in regularization theory~\cite{EHN96}.
A standard technique to compute a stable solution for an ill-conditioned
system is to use a truncated singular value decomposition
(TSVD)~\cite{EHN96}. 
This means in our case we compute a regularized
pseudo-inverse $\Rnpe = V_n \Sigma_n^{+,\thresh} U_n^{\ast}$  where
\begin{equation}
\Sigma^{+,\thresh} = \diag(\{d_{k}^{+}\})\,, \quad 
d_{k}^{+} = 
\begin{cases}
1/\lambda_k & \text{if} \,\,\lambda_k \ge \thresh, \\
0 & \text{otherwise.}
\end{cases}
\label{pinvtrunc}
\end{equation}
In~\cite{Har98} it is shown that for each $n$ we can choose
an appropriate truncation level $\thresh$ such that the regularized 
inverses $\Rnpe$ converge strongly to $\Rp$ for $\ntoinf$ and 
consequently $\underset{\ntoinf}{\lim} \|f - \fn\| = 0$, where 
\begin{equation}
\fn(t) = \sum_{j=-n}^{n} c_{j}^{(n,\thresh)} \sinco(t-t_j) \notag
\end{equation}
with
\begin{equation}
c^{(n,\thresh)} = \Rnpe \bn\,. \notag
\end{equation}
The optimal truncation level $\thresh$ depends on the dimension $n$, the
sampling geometry, and the noise level. Thus it is not known a priori
and has in principle to be determined for each $n$ independently. 

Since $\thresh$ is of vital importance for the quality
of the reconstruction, but no theoretical explanations
for the choice of $\thresh$ are given in the sampling literature,  
we briefly discuss this issue.
For this purpose we need some results from regularization theory.

\subsection{Estimation of regularization parameter} \label{ss:est}

Let $Ax = y^{\delta}$ be given where $A$ is ill-conditioned or singular
and $y^{\delta}$ is a perturbed right-hand side with 
$\|y - y^{\delta}\| \le \delta \|y\|$.
Since in our sampling problem the matrix under consideration is symmetric, 
we assume for convenience that $A$ is symmetric.
From a numerical point of view ill-conditioned systems behave like
singular systems and additional information is needed to obtain
a satisfactory solution to $Ax=y$. This information is usually stated in terms
of ``smoothness'' of the solution $x$. A standard approach to qualitatively 
describe smoothness of $x$ is to require that $x$ can be represented in 
the form $x=Sz$ with some vector $z$ of reasonable norm, and a ``smoothing'' 
matrix $S$, cf.~\cite{EHN96,Neu98}. Often it is useful to construct 
$S$ directly from $A$ by setting
\begin{equation}
S = A^p\,,   \qquad p \in \Nst_0 \,.
\label{smoothness}
\end{equation}
Usually, $p$ is assumed to be fixed, typically at $p=1$ or $p=2$.

We compute a regularized solution to $Ax=y^{\delta}$ via a truncated
SVD and want to determine the optimal regularization
parameter (i.e., truncation level) $\tau$.

Under the assumption that 
\begin{equation}
x = Sz \, , \quad  \|Ax-y^{\delta}\| \le \Delta \|z\|
\label{}
\end{equation}
it follows from Theorem~4.1 in ~\cite{Neu98} that the optimal
regularization parameter $\tau$ for the TSVD is 
\begin{equation}
\hat{\tau}=\left(\frac{\gamma_1 \delta}{\gamma_2 p}\right)^{\frac{1}{p+1}}\,,
\label{optreg}
\end{equation}
where $\gamma_1=\gamma_2 =1$ (see Section~6 in~\cite{Neu98}).

However $z$ and $\Delta$ are in general not known. Using
$\|Ax-y^{\delta}\| \le \delta \|y\|$ and $\|y\|=\|Ax\|=\|A Sz\|=\|A^{p+1} z\|$
we obtain $ \|y\| \le \|A\|^{p+1}\|z\|$. Furthermore, setting
$\delta \|y\| = \Delta \|z\|$ implies 
\begin{equation}
\Delta \le \delta \|A\|^{p+1}\,.
\label{Delta}
\end{equation}
Hence combining~\eqref{optreg} and~\eqref{Delta} we get
\begin{equation}
\hat{\tau} \le \left( \frac{\delta \|A\|^{p+1}}{p}\right)^{\frac{1}{p+1}} 
= \|A\| \left(\frac{\delta}{p}\right)^{\frac{1}{p+1}}\,.
\label{tauest}
\end{equation}

Applying these results to solving $\Rn \cn = \bn$ via TSVD as described
in the previous section, we get
\begin{equation}
\hat{\tau} \le \|\Rn\| \left(\frac{\delta}{p}\right)^{\frac{1}{p+1}}
\le \|R\| \left(\frac{\delta}{p}\right)^{\frac{1}{p+1}}
= B \left(\frac{\delta}{p}\right)^{\frac{1}{p+1}}\,,
\label{threshopt}
\end{equation}
where $B$ is the upper frame bound. Fortunately estimates for the upper
frame bound are much easier to obtain than estimates for the lower
frame bound.

Thus using the standard setting $p=1$ or $p=2$ a good choice for the
regularization parameter $\tau$ is
\begin{equation}
\thresh \subseteq [B(\delta/2)^{1/3},B(\delta)^{1/2}]\,.
\label{thresh}
\end{equation}
Extensive numerical simulations confirm this choice, see also
Section~\ref{s:applications}.

For instance for the reconstruction problem
of Example~1 with noise-free data and machine
precision $\eps = \delta = 10^{-16}$, formula~\eqref{thresh} implies 
$\thresh \subseteq [10^{-6},10^{-8}]$. This coincides very well
with numerical experiments.

If the noise level $\delta$ is not known, it has to be estimated.
This difficult problem will not be discussed here.
The reader is referred to~\cite{Neu98} for more details.

Although we have arrived now at an implementable algorithm for
the nonuniform sampling problem, the disadvantages of the approach described
in the previous section are obvious.
In general the matrix $\Rn$ does not have any particular structure, thus 
the computational costs for the singular value decomposition are $\ord(n^3)$ 
which is prohibitive large in many applications. 
It is definitely not a good approach to transform a well-posed
infinite-dimensional problem into an ill-posed finite-dimensional problem
for which a stable solution can only be computed by using a ``heavy
regularization machinery''.

The methods in~\cite{Yen56,YT67,Win92,San94,BH90} coincide with or are
essentially equivalent to the truncated frame approach, therefore they
suffer from the same instability problems and the same
numerical inefficiency.

\subsection{CG and regularization of the truncated frame method}
\label{ss:cgtrunc}

As mentioned above one way to stabilize the solution
of $\Rn \cn = \bn$ is a truncated singular value decomposition, where
the truncation level serves as regularization parameter.
For large $n$ the costs of the singular value decomposition become 
prohibitive for practical purposes. 

We propose the conjugate gradient method~\cite{GL96} to solve $\Rn \cn = \bn$. 
It is in general much more efficient than a TSVD (or Tikhonov
regularization as suggested in~\cite{Win92}), and at the same time it 
can also be used as a regularization method.

The standard error analysis for CG cannot be used in our case, since
the matrix is ill-conditioned. Rather we have to resort to the
error analysis developed in~\cite{NP84,Han95}.

When solving a linear system $Ax=y$ by CG for noisy data
$y^{\delta}$ following happens. The iterates $x_k$ of CG
may diverge for $k \rightarrow \infty$, however the error propagation
remains limited in the beginning of the iteration. The quality of the
approximation therefore depends on how many iterative steps can be
performed until the iterates turn to diverge. The idea is now to stop
the iteration at about the point where divergence sets in.
In other words the iterations count is the regularization parameter which
remains to be controlled by an appropriate stopping rule~\cite{Nem86,Han95}.

In our case assume $\|\bnd-\bn\| \le \delta \|\bn\|$, where
$b_j^{(n,\delta)}$ denotes a noisy sample. 
We terminate the CG iterations when the iterates $(\cnd)_k$ satisfy for the 
first time~\cite{Han95}
\begin{equation}
\|\bn\ - (\cnd)_k\| \le \tau \delta \|\bn\| 
\label{stopcg}
\end{equation}
for some fixed $\tau >1$.

It should be noted that one can construct ``academic'' examples where this 
stopping rule does not prevent CG from diverging, see~\cite{Han95},
``most of the time'' however it gives satisfactory results.
We refer the reader to~\cite{Nem86,Han95} for a detailed discussion of
various stopping criteria.

There is a variety of reasons, besides the ones we have already mentioned,
that make the conjugate gradient method and
the nonuniform sampling problem a ``perfect couple''.
See Sections~\ref{ss:trigpol},~\ref{ss:ml}, and~\ref{ss:regul}
for more details.

By combining the truncated frame approach with the conjugate gradient method
(with appropriate stopping rule) we finally arrive at a reconstruction
method that is of some practical relevance. However the only existing
method at the moment that can handle large scale reconstruction problems
seems to be the one proposed in the next section.

\section{Trigonometric polynomials and efficient signal reconstruction} 
\label{ss:trigpol}

In the previous section we have seen that the naive finite-dimensional
approach via truncated frames is not satisfactory, it already
leads to severe stability problems in the ideal case of regular
oversampling. In this section we propose a different finite-dimensional
model, which resembles much better the structural properties of the 
sampling problem, as can be seen below.

The idea is simple. In practice only a finite number of samples 
$\{f(t_j)\}_{j=1}^{r}$ is given, where without loss of generality we assume 
$-M \le t_1 < \dots < t_r \le M$ (otherwise we can always re-normalize the 
data).  Since no data of $f$ are available from outside this region
we focus on a local approximation of $f$ on $[-M,M]$.
We extend the sampling set periodically across the boundaries, 
and identify this interval with the (properly normalized)
torus $\Tst$. To avoid technical problems at the boundaries in the sequel 
we will choose the interval somewhat larger and consider either 
$[-M-1/2,M+1/2]$ or $[-N,N]$ with $N= M+\frac{M}{r-1}$. For
theoretical considerations the choice $[-M-1/2,M+1/2]$ is more convenient.

Since the dual group of the torus $\Tst$ is $\Zst$, periodic 
band-limited functions on $\Tst$ reduce to trigonometric polynomials
(of course technically $f$ does then no longer belong to $\BO$ since
it is no longer in $\LtR$).
This suggests to use trigonometric polynomials as a realistic
finite-dimensional model for a numerical solution of the nonuniform
sampling problem. We consider the space $\PM$ of
trigonometric polynomials of degree $M$ of the form
\begin{equation}
p(t) = (2M+1)^{-1} \sum_{k=-M}^{M} a_{k} e^{2\pi i kt/(2M+1)}\,.
\label{pm}
\end{equation}
The norm of $p \in \PM$ is 
$$\|p\|^2 =\int \limits_{-N}^{N} |p(t)|^2\, dt =\sum_{k=-M}^{M} |a_k|^2 \,.$$
Since the distributional Fourier transform of $p$ is 
$\hat{p} = (2M+1)^{-1} \sum_{k=-M}^{M} a_k \delta_{k/(2M+1)}$ we have
$\supp \hat{p} \subseteq \{ k/(2M+1) , |k| \le M\} \subseteq [-1/2, 1/2]$.
Hence $\PM$ is indeed a natural finite-dimensional model for $\BO$.


In general the $f(t_j)$ are not the samples of a trigonometric
polynomial in $\PM$, moreover the samples are usually perturbed by noise,
hence we may not find a $p \in \PM$ such that $p(t_j)=b_j = f(t_j)$. We 
therefore consider the least squares problem
\begin{equation}
\underset{p \in \PM }{ \min } \sum_{j=1}^{r} |p(t_j) -b_j|^2 w_j\,.
\label{LSP}
\end{equation}
Here the $w_j >0$ are user-defined weights, which can be chosen for instance 
to compensate for irregularities in the sampling geometry~\cite{FGS95}.

By increasing $M$ so that $ r \le 2M+1$ we can certainly find
a trigonometric polynomial that interpolates the given data exactly.
However in the presence of noise, such a solution is usually rough and
highly oscillating and may poorly resemble the original signal.
We will discuss the question of the optimal choice of $M$ if the original
bandwidth is not known and in presence of noisy data in Section~\ref{ss:regul}.

The following theorem provides an efficient numerical 
reconstruction algorithm. It is also the key for the analysis of the
relation between the finite-dimensional approximation in $\PM$ and
the solution of the original infinite-dimensional sampling problem
in $\BO$.

\begin{theorem}[{\bf and Algorithm}] \cite{Gro93a,FGS95}
\label{th:act}
Given the sampling points $-M \le t_1 < \dots, t_{r} \le M$, samples
$\{b_j\}_{j=1}^r$, positive weights $\{w_j\}_{j=1}^{r}$
with $2M+1 \le r$. \\
Step 1: Compute the $(2M+1)\times (2M+1)$ Toeplitz matrix $T_M$ 
with entries
\begin{equation}
(T_M)_{k,l}=\frac{1}{2M+1}\sum_{j=1}^{r} w_j e^{-2\pi i (k-l) t_j/(2M+1)}
    \qquad \mbox{for $|k|,|l| \le M$}
\label{toepmat}
\end{equation}
and $y_M \in \Cst^{(2M+1)}$ by
\begin{equation}
\label{rightside}
(y_M)_k=\frac{1}{\sqrt{2M+1}} \sum_{j=1}^{r} b_j w_j e^{-2\pi i k t_j/(2M+1)}
    \qquad \mbox{for $|k| \le M$} \,.
\end{equation}
\noindent Step 2: Solve the system
\begin{equation}
\label{toepsys}
T_M a_M = y_M \,.
\end{equation}
\noindent Step 3: Then the polynomial $\plsp \in \PM$ that
solves~\eqref{LSP} is given by
\begin{equation}
\plsp(t)=\frac{1}{\sqrt{2M+1}} \sum_{k=-M}^M (a_M)_k e^{2 \pi i kt/(2M+1)} \,.
\label{lsppol}
\end{equation}
\end{theorem}

\noindent
{\bf Numerical Implementation of Theorem/Algorithm~\ref{th:act}:} \\
Step 1: The entries of $T_M$ and $y_M$ of equations~\eqref{toepmat} 
and~\eqref{rightside}
can be computed in $\ord(M \log M + r \log(1/\eps))$ operations (where
$\eps$ is the required accuracy) using Beylkin's unequally spaced FFT 
algorithm~\cite{Bey95}.\\
Step 2: We solve $T_M a_M = y_M$ by the conjugate gradient (CG)
algorithm~\cite{GL96}. The matrix-vector multiplication in each iteration
of CG can be carried out in $\ord (M \log M)$ operations via 
FFT~\cite{CN96}. Thus the solution of~\eqref{toepsys} takes $\ord (k M \log M)$
operations, where $k$ is the number of iterations.\\
Step 3: Usually the signal is reconstructed on regularly space nodes
$\{u_i\}_{i=1}^{N}$. In this case $p_M(u_i)$ in~\eqref{lsppol} can be
computed by FFT. For non-uniformly spaced nodes $u_i$ we can again resort
to Beylkin's USFFT algorithm.  
 
\medskip
There exists a large number of fast algorithms for the solution of Toeplitz
systems. Probably the most efficient algorithm in our case is CG.
We have already mentioned that the Toeplitz system~\eqref{toepsys}
can be solved in $\ord (kM \log M)$ via CG. The number of iterations $k$
depends essentially on the clustering of the eigenvalues of $T_M$, 
cf.~\cite{CN96}. It follows from equation~\eqref{kron} below and perturbation 
theory~\cite{Chr96a} that, if the sampling points stem from a perturbed regular 
sampling set, the eigenvalues of $T_M$ will be clustered around $\beta$,
where $\beta$ is the oversampling rate. In such cases we can expect
a very fast rate of convergence. 
The simple frame iteration~\cite{Mar93a,Ben92} is not able to take 
advantage of such a situation.

\medskip

For the analysis of the relation between the solution $\plsp$ of
Theorem~\ref{th:act} and the solution $f$ of the original infinite-dimensional 
problem we follow Gr\"ochenig~\cite{Gro99}. Assume that 
the samples $\{f(t_j)\}_{j \in \Zst}$ of $f \in \BO$ are given. For the 
finite-dimensional approximation we consider
only those samples $f(t_j)$ for which $t_j$ is contained in the interval 
$[-M-\frac{1}{2}, M+\frac{1}{2}]$ and compute the least squares
approximation $\plsp$ with degree $M$ and period $2M+1$ as
in Theorem~\ref{th:act}. It is shown in~\cite{Gro99} that if 
$\sigma (T_M) \subseteq [\alpha, \beta]$ for all $M$ with $\alpha >0$ then
\begin{equation}
\underset{M \toinf}{\lim} \int \limits_{[-M, M]} |f(t) - \plsp(t)|^2 \, dt =0 ,
\label{trigconv}
\end{equation}
and also $\lim \plsp(t) = f(t)$ uniformly on compact sets.

Under the Nyquist condition $\sup(t_{j+1}-tj) :=\gamma < 1$ and using
weights $w_j = (t_{j+1}-t_{j-1})/2$ Gr\"ochenig has shown that
\begin{equation}
\sigma (T_M) \subseteq [(1-\gamma)^2, 6]\,,
\label{condest}
\end{equation}
independently of $M$, see~\cite{Gro99}. These results validate the
usage of trigonometric polynomials as finite-dimensional model
for nonuniform sampling.

\medskip

\noindent
{\bf Example 1 -- reconsidered:}
Recall that in Example~1 of Section~\ref{ss:truncated} we have considered the
reconstruction of a regularly oversampled signal $f \in \BO$.
What does the reconstruction method 
of Theorem~\ref{th:act} yield in this case? Let us check the entries of 
the matrix $T_M$ when we take only those samples in the interval $[-n,n]$. 
The period of the polynomial becomes $2N$ with $N=n+\frac{n}{r-1}$ 
where $r$ is the number of given samples. Then
\begin{equation}
(T_M)_{k,l} = \frac{1}{2N} \sum_{j=1}^{r} e^{2\pi i (k-l)t_j/(2N)}
 = \sum_{j=-nm}^{nm} e^{2\pi i (k-l) \frac{j}{2nm+1}} = m\delta_{k,l}
\label{kron}
\end{equation}
for $k,l = -M,\dots,M$, 
where $\delta_{k,l}$ is Kronecker's symbol with the usual
meaning $\delta_{k,l}=1$ if $k=l$ and $0$ else. Hence we get
$$T_M = m I\,,$$
where $I$ is the identity matrix on $\Cst^{2M+1}$, thus $T_M$ 
resembles the structure of the infinite-dimensional frame operator $S$ 
in this case (including exact approximation of the frame bounds).
Recall that the truncated frame approach leads to an ``artificial'' 
ill-posed problem even in such a simple situation.

The advantages of the trigonometric polynomial approach compared to
the truncated frame approach are manifold.
In the one case we have to deal with an ill-posed
problem which has no specific structure, hence its solution is
numerically very expensive. In the other case we have to solve a
problem with rich mathematical structure, whose stability depends only on 
the sampling density, a situation that resembles the original 
infinite-dimensional sampling problem. 
 
\medskip

In principle the coefficients $\alsp=\{\alspk\}_{k=-M}^{M}$ of the polynomial 
$\plsp$ that minimizes~\eqref{LSP} could also be computed by directly
solving the Vandermonde type system
\begin{equation}
WV \alsp = W b\,,
\label{vandermonde}
\end{equation}
where $V_{j,k}=\frac{1}{\sqrt{2M+1}} e^{-2 \pi i k t_j/(2M+1)}$ for
$j=1,\dots, r,\, k=-M,\dots, M$ and $W$ is a diagonal matrix with
entries $W_{j,j}= \sqrt{w_j}$, cf.~\cite{RAG91}. Several algorithms are 
known for a relatively efficient solution of Vandermonde 
systems~\cite{BP70,RAG91}. 
However this is one of the rare cases, where, instead of directly 
solving~\eqref{vandermonde}, it is advisable to explicitly establish 
the system of normal equations 
\begin{equation}
T_M a_M = y_M\,,
\label{normal}
\end{equation}
where $T=V^{\ast} W^2 V$ and $y = V^{\ast} W^2 b$.

The advantages of considering the system $T_M a_M = y_M$ instead of
the Vandermonde system~\eqref{vandermonde} are manifold:
\begin{itemize}
\vspace*{-1mm}
\setlength{\itemsep}{-0.5ex}
\setlength{\parsep}{-0.5ex}
\item The matrix $T_M$ plays a key role in the analysis of the relation of 
the solution of~\eqref{LSP} and the solution of the infinite-dimensional 
sampling problem~\eqref{recon1}, see~\eqref{trigconv} and ~\eqref{condest} 
above. 
\item $T_M$ is of size $(2M+1) \times (2M+1)$, independently of the number of 
sampling points. Moreover, since 
$(T_M)_{k,l}=\sum_{j=1}^{r} w_j e^{2\pi i (k-l) t_j}$,
it is of Toeplitz type. 
These facts give rise to fast and robust reconstruction algorithms.
\item The resulting reconstruction algorithms can be easily generalized to 
higher dimensions, see Section~\ref{ss:multi}. Such a generalization to 
higher dimensions seems not 
to be straightforward for fast solvers of Vandermonde systems such 
as the algorithm proposed in~\cite{RAG91}.
\end{itemize}

\medskip

We point out that other finite-dimensional approaches are proposed
in~\cite{FLS98,CC98}. These approaches may provide interesting
alternatives in the few cases where the algorithm outlined in 
Section~\ref{ss:trigpol} does not lead to good results.
These cases occur when only a few samples of the signal $f$ are given
in an interval $[a,b]$ say, and at the same time we have 
$|f(a) - f(b)| \gg 0$ and $|f'(a) - f'(b)| \gg 0$, i.e., if
$f$ is ``strongly non-periodic'' on $[a,b]$.
However the computational complexity of the methods in~\cite{FLS98,CC98} is 
significantly larger.

\subsection{Multi-dimensional nonuniform sampling} \label{ss:multi}

The approach presented above can be easily generalized to higher 
dimensions by a diligent book-keeping of the notation.
We consider the space of $d$-dimensional trigonometric polynomials $\PMd$
as finite-dimensional model for $\BO^d$. For given samples $f(t_j)$ of 
$f \in \BO^d$, where $t_j \in \Rdst$, we compute the least squares 
approximation $\plsp$ similar to Theorem~\ref{th:act} by solving the 
corresponding system of equations $T_M a_M = y_M$. 

In 2-D for instance the matrix $T_M$ 
becomes a block Toeplitz matrix with Toeplitz blocks~\cite{Str97}.
For a fast computation of the entries of $T$ we can again make use
of Beylkin's USFFT algorithm~\cite{Bey95}. And similar to 1-D, 
multiplication of a vector by $T_M$ can be carried out by 2-D FFT. 

Also the relation between the finite-dimensional approximation in
$\PMd$ and
the infinite-dimensional solution in $\BO^{d}$ is similar as in 1-D. 
The only mathematical difficulty is to give conditions under which
the matrix $T_M$ is invertible. Since the fundamental theorem of
algebra does not hold in dimensions larger than one, the condition
$(2M+1)^d \le r$ is necessary but no longer sufficient for
the invertibility of $T_M$. Sufficient conditions for the invertibility,
depending on the sampling density, are presented in~\cite{Gro99a}. 

\section{Bandwidth estimation and regularization} \label{s:numeric}

In this section we discuss several numerical aspects of nonuniform sampling
that are very important from a practical viewpoint, however only few
answers to these problems can be found in the literature.

\subsection{A multilevel signal reconstruction algorithm} 
\label{ss:ml}

In almost all theoretical results and numerical algorithms for reconstructing
a band-limited signal from nonuniform samples it is assumed that the
bandwidth is known a priori. This information however is often not 
available in practice. 

A good choice of the bandwidth for the reconstruction algorithm becomes 
crucial in case of noisy data.
It is intuitively clear that choosing a too large bandwidth leads to 
over-fit of the noise in the data, while a too small bandwidth yields
a smooth solution but also to under-fit of the data. And of course we
want to avoid the determination of the ``correct'' $\Omega$ by 
trial-and-error methods. Hence the problem is to design a method
that can reconstruct a signal from non-uniformly spaced, noisy samples
without requiring a priori information about the bandwidth of the signal.

The multilevel approach derived in~\cite{SS97} provides an answer
to this problem. The approach applies to an infinite-dimensional
as well as to a finite-dimensional setting. We describe the method
directly for the trigonometric polynomial model, where the determination
of the bandwidth $\Omega$ translates into the determination of the polynomial
degree $M$ of the reconstruction. The idea of the multilevel algorithm
is as follows. 

Let the noisy samples $\{b^{\delta}_j\}_{j=1}^{r}=\{f^{\delta}(t_j)\}_{j=1}^r$ 
of $f \in \BO$ be given with  
$\sum_{j=1}^{r}|f(t_j)-b^{\delta}(t_j)|^2 \le \delta^2 \|b^{\delta}\|^2$ and
let $Q_M$ denote the orthogonal projection from $\BO$ into $\PM$.
We start with initial degree $M=1$ and run Algorithm~\ref{th:act}
until the iterates $p_{0,k}$ satisfy for the first time the {\em inner}
stopping criterion
\begin{equation}
\sum_{j=1}^{r}|p_{1,k}(t_j) - b^{\delta}_j|^2 \le 
2 \tau (\delta \|b^{\delta}\| + \|Q_0 f - f\|)\|b^{\delta}\|
\notag
\end{equation}
for some fixed $\tau >1$.
Denote this approximation (at iteration $k_*$) by $p_{1,k_*}$.
If $p_{1,k_*}$ satisfies the {\em outer} stopping criterion
\begin{equation}
\sum_{j=1}^{r}|p_{1,k}(t_j) - b^{\delta}_j|^2 \le 
2 \tau \delta \|b^{\delta}\|^2
\label{stopout}
\end{equation}
we take $p_{1,k_*}$ as final approximation. Otherwise we
proceed to the next level $M=2$ and run Algorithm~\ref{th:act}
again, using $p_{1,k_*}$ as initial approximation
by setting $p_{2,0} =p_{1,k_*}$. 

At level $M=N$ the inner level-dependent stopping criterion becomes
\begin{equation}
\sum_{j=1}^{r}|p_{N,k}(t_j) - b^{\delta}_j|^2 \le 
2 \tau (\delta \|b^{\delta}\| + \|Q_N f - f\|)\|b^{\delta}\|,
\label{stopin}
\end{equation}
while the outer stopping criterion does not change since it
is level-independent. 

Stopping rule~\eqref{stopin} guarantees that the iterates of CG do not diverge.
It also ensures that CG does not iterate too long at a certain level, since
if $M$ is too small further iterations at this level will not lead
to a significant improvement. Therefore we switch to the next level.
The outer stopping criterion~\eqref{stopout} controls over-fit and
under-fit of the data, since in presence of noisy data is does not make sense
to ask for a solution $p_M$ that satisfies
$\sum_{j=1}^{r}|p_{M}(t_j) - b^{\delta}_j|^2=0$.

Since the original signal $f$ is not known, the expression
$\|f - Q_N f\|$ in~\eqref{stopin} cannot be computed.  In~\cite{SS97} the 
reader can find an approach to estimate $\|f - Q_N f\|$ recursively.

\subsection{Solution of ill-conditioned sampling problems} 
\label{ss:regul}

A variety of conditions on the sampling points $\{t_j\}_{j \in \Zst}$ are 
known under which the set $\sincframe$ is a frame for $\BO$, which in turn 
implies (at least theoretically) perfect reconstruction of a signal $f$ from 
its samples $f(t_j)$. This does however
not guarantee a stable reconstruction from a numerical viewpoint, since the 
ratio of the frame bounds $B/A$ can still be extremely large and therefore 
the frame operator $S$ can be ill-conditioned. This may happen for
instance if $\gamma$ in~\eqref{condest} goes to 1, in which case 
$\cond(T)$ may become large. The sampling problem may also become numerically 
unstable or even ill-posed, if the sampling set has large gaps, which is very 
common in astronomy and geophysics. Note that in this case the instability of 
the system $T_M a_M = y_M$ does {\em not} result from an inadequate 
discretization of the infinite-dimensional problem.

There exists a large number of (circulant) Toeplitz preconditioners that 
could be applied to the system $T_M a_M = y_M$, however it turns out that 
they do not improve the stability of the problem in this case. The reason 
lies in the distribution of the eigenvalues of $T_M$, as we will see below.

Following~\cite{Tyr96}, we call two sequences of real numbers
$\{\lambda^{(n)}\}_{k=1}^{n}$ and $\{\nu^{(n)}\}_{k=1}^{n}$ 
{\em equally distributed}, if
\begin{equation}
\underset{\ntoinf}{\lim} \frac{1}{n} \sum_{k=1}^{n}
[F(\lambda^{(n)}_{k}) - F(\nu^{(n)}_{k}) ] = 0
\label{defdist}
\end{equation}
for any continuous function $F$ with compact support\footnote{In H.Weyl's 
definition $\lambda^{(n)}_{k}$ and $\nu^{(n)}_{k}$ are required to belong to
a common interval.}.

Let $C$ be a $(n \times n)$ circulant matrix with first column 
$(c_0,\dots,c_{n-1})$, we write $C = \circ (c_0,\dots,c_{n-1})$. The 
eigenvalues of $C$ are distributed as 
$\lambda_k = \frac{1}{\sqrt{n}}\sum_{l=0}^{n-1} c_l e^{2\pi i kl/n}$.
Observe that the Toeplitz matrix $A_n$ with first
column $(a_0,a_1,\dots, a_n)$ can be embedded in the circulant
matrix 
\begin{equation}
C_n =\circ (a_0,a_1,\dots, a_n, \bar{a_n},\dots, \bar{a_1})\,.
\label{circembed}
\end{equation}
Thms~4.1 and~4.2 in~\cite{Tyr96} state that the eigenvalues of $A_n$ and
$C_n$ are equally distributed as $f(x)$ where
\begin{equation}
f(x) = \sum_{k=-\infty}^{\infty} a_k e^{2 \pi i kx}\,.
\label{fcirc}
\end{equation}
The partial sum of the series~\eqref{fcirc} is
\begin{equation}
f_n(x) = \sum_{k=-n}^{n} a_k e^{2 \pi i kx}\,.
\label{fcircm}
\end{equation}

To understand the clustering behavior of the eigenvalues of $T_M$ in
case of sampling sets with large gaps, we consider a sampling set in
$[-M,M)$, that consists of one large block of samples and one large gap, i.e., 
$t_j = \frac{j}{Lm}$ for $j=-mM,\dots mM$ for $m,L \in \Nst$.
(Recall that we identify the interval with the torus). Then the entries 
$z_k$ of the Toeplitz matrix $T_M$ of~\eqref{toepmat} (with $w_j=1$) are 
$$z_k=\frac{1}{2M+1}\sum_{j=-mM}^{mM} e^{-2\pi i k \frac{j}{Lm}/(2M+1)},
\quad k=0,\dots,2M\,.$$ 
To investigate the clustering behavior of the eigenvalues of $T_M$ for 
$M \toinf$, we embed $T_M$ in a circulant matrix $C_M$ as 
in~\eqref{circembed}. Then~\eqref{fcircm} becomes
\begin{equation}
f_{mM}(x) = \frac{1}{Lm(2M+1)}\sum_{l=-mM}^{mM} \sum_{j=-mM}^{mM} 
e^{2 \pi il [k/(4M+1) - j/((2M+1)mL)]}
\end{equation}
whence $f_{mM} \rightarrow {\bf 1}_{[-1/(2L),1/(2L)]}$ for $M \toinf$, where
${\bf 1}_{[-a,a]}(x) = 1$, if $-a  < x < a$ and 0 else.

Thus the eigenvalues of $T_M$ are asymptotically clustered around zero and one.
For general nonuniform sampling sets with large gaps the clustering at 1
will disappear, but of course the spectral cluster at 0 will remain.
In this case it is known that the preconditioned problem will still have a 
spectral cluster at the origin~\cite{YC93} and preconditioning will not
be efficient.


\medskip

Fortunately there are other possibilities to obtain a stabilized
solution of $T_M a_M = y_M$. The condition number of $T_M$ essentially depends
on the ratio of the maximal gap in the sampling set to
the Nyquist rate, which in turn depends on the bandwidth of the signal.
We can improve the stability of the system by adapting the degree $M$
of the approximation accordingly. Thus the parameter $M$ serves as a
regularization parameter that balances stability and accuracy of
the solution. This technique can be seen as a specific realization of 
{\em regularization by projection}, see Chapter~3 in~\cite{EHN96}.
In addition, as described in Section~\ref{ss:regul}, we can utilize CG 
as regularization method for the solution of the Toeplitz system in order 
to balance approximation error and propagated error. The multilevel method 
introduced in Section~\ref{ss:ml} combines both features. By optimizing the 
level (bandwidth) and the number of iterations in each level it provides
an efficient and robust regularization technique for ill-conditioned
sampling problems. See Section~\ref{s:applications} for numerical
examples.

\section{Applications} \label{s:applications}

We present two numerical examples to demonstrate the performance of the
described methods.
The first one concerns a 1-D reconstruction problem arising in 
spectroscopy. In the second example we approximate
the Earth's magnetic field from noisy scattered data.

\subsection{An example from spectroscopy} \label{ss:spectro}

The original spectroscopy signal $f$ is known at 1024 regularly spaced 
points $t_j$. This discrete sampling sequence will play the role of the 
original continuous signal.
To simulate the situation of a typical experiment in spectroscopy we
consider only 107 randomly chosen sampling values of the given sampling set. 
Furthermore we add noise to the samples with noise level (normalized by
division by $\sum_{k=1}^{1024}|f(t_j)|^2$) of $\delta=0.1$. 
Since the samples are contaminated by noise, we cannot expect
to recover the (discrete) signal $f$ completely. The bandwidth is 
approximately $\Omega =5$ which translates into a polynomial degree of
$M \approx 30$. Note that in general $\Omega$ and (hence $M$) may not be 
available. We will also consider this situation, but in the first
experiments we assume that we know $\Omega$. The error between the
original signal $f$ and an approximation $f_n$ is measured by
computing $\|f- f_n\|^2/\|f\|^2$.

\begin{figure}
\begin{center}
\subfigure[Truncated frame method with TSVD,
error=0.0944.]{
\epsfig{file=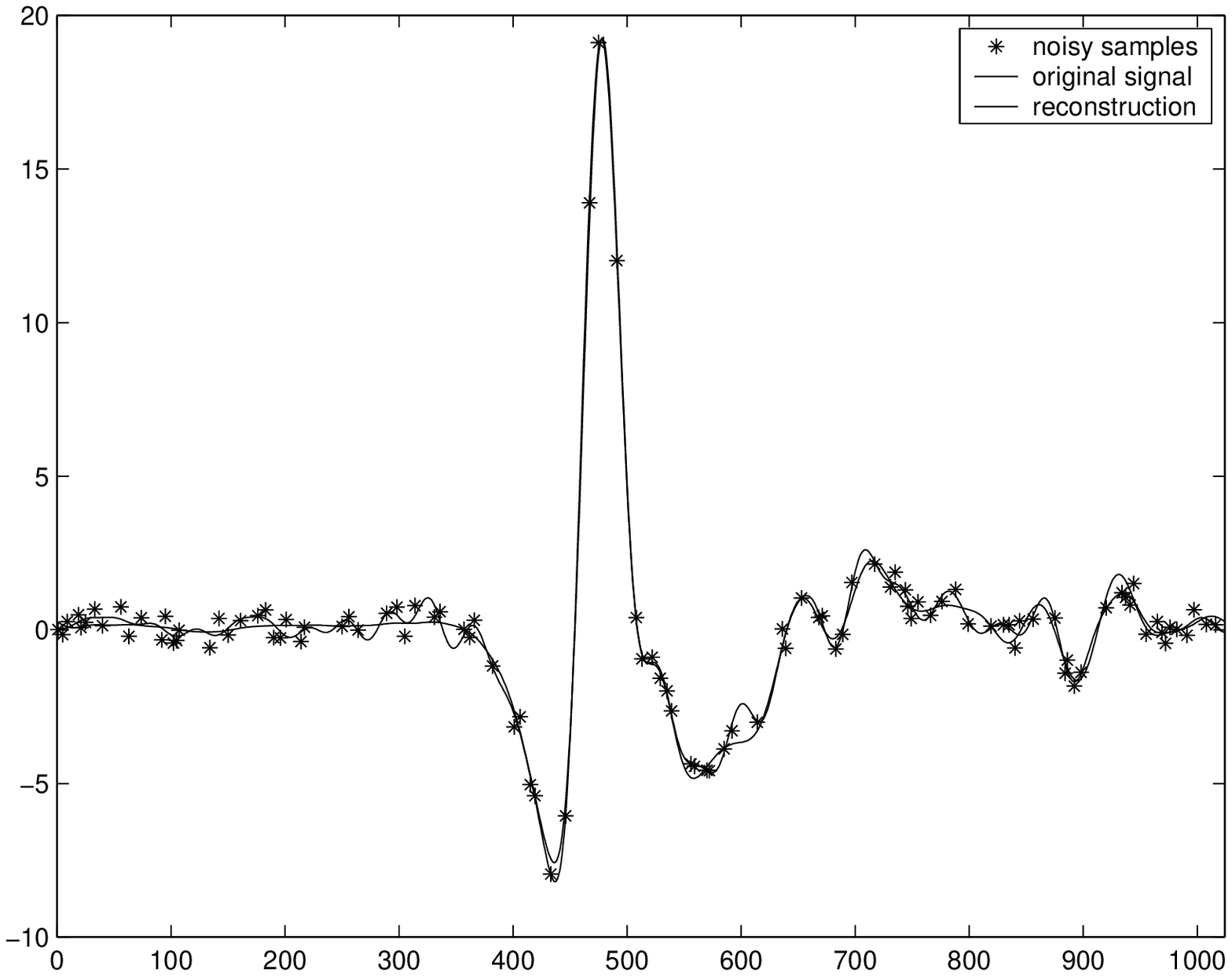,width=55mm,height=40mm}}
\subfigure[Truncated frame method with CG, 
error=0.1097.]{
\epsfig{file=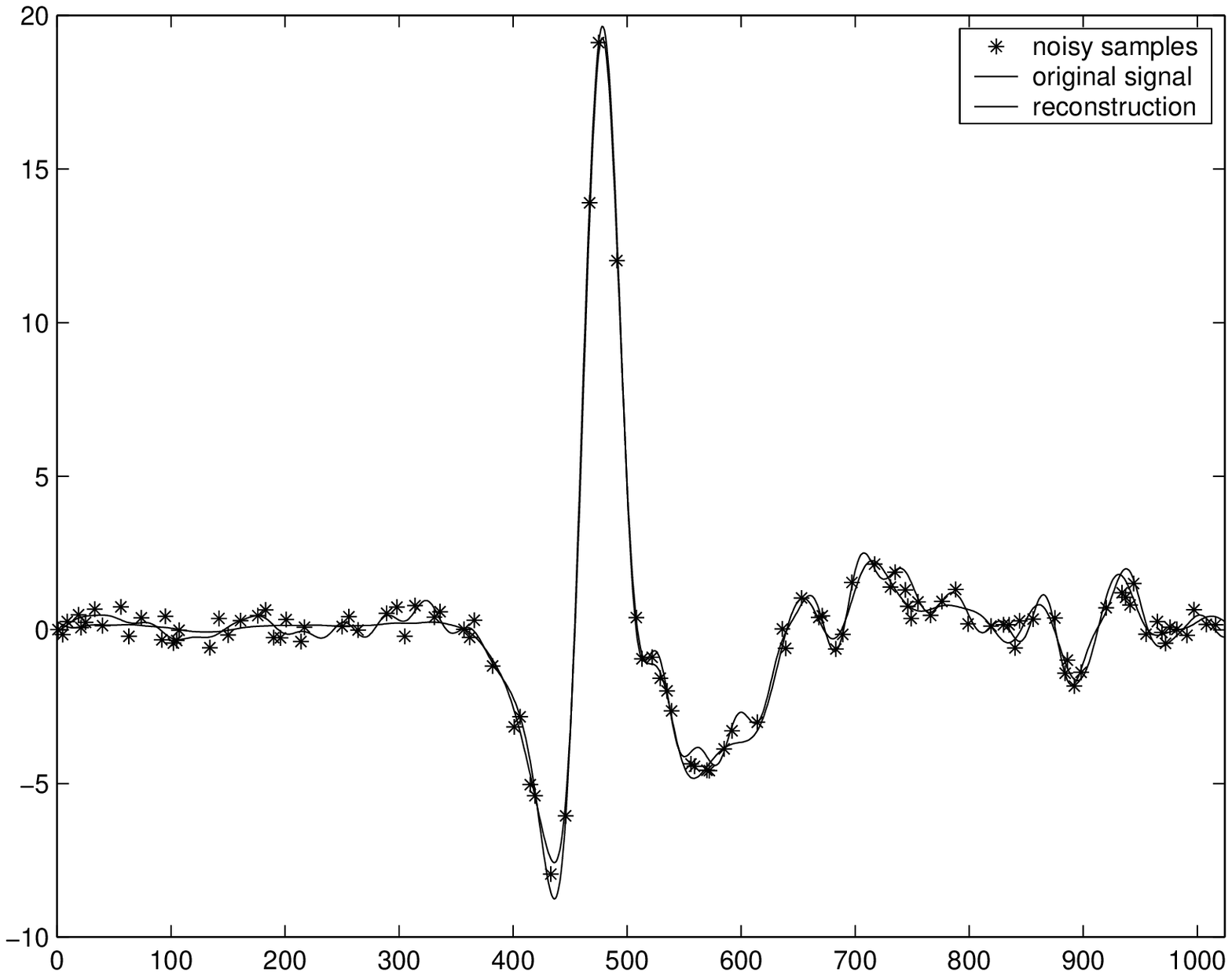,width=55mm,height=40mm}}\\
\subfigure[Algorithm~\ref{th:act} with ``correct'' bandwidth,
error=0.0876]{
\epsfig{file=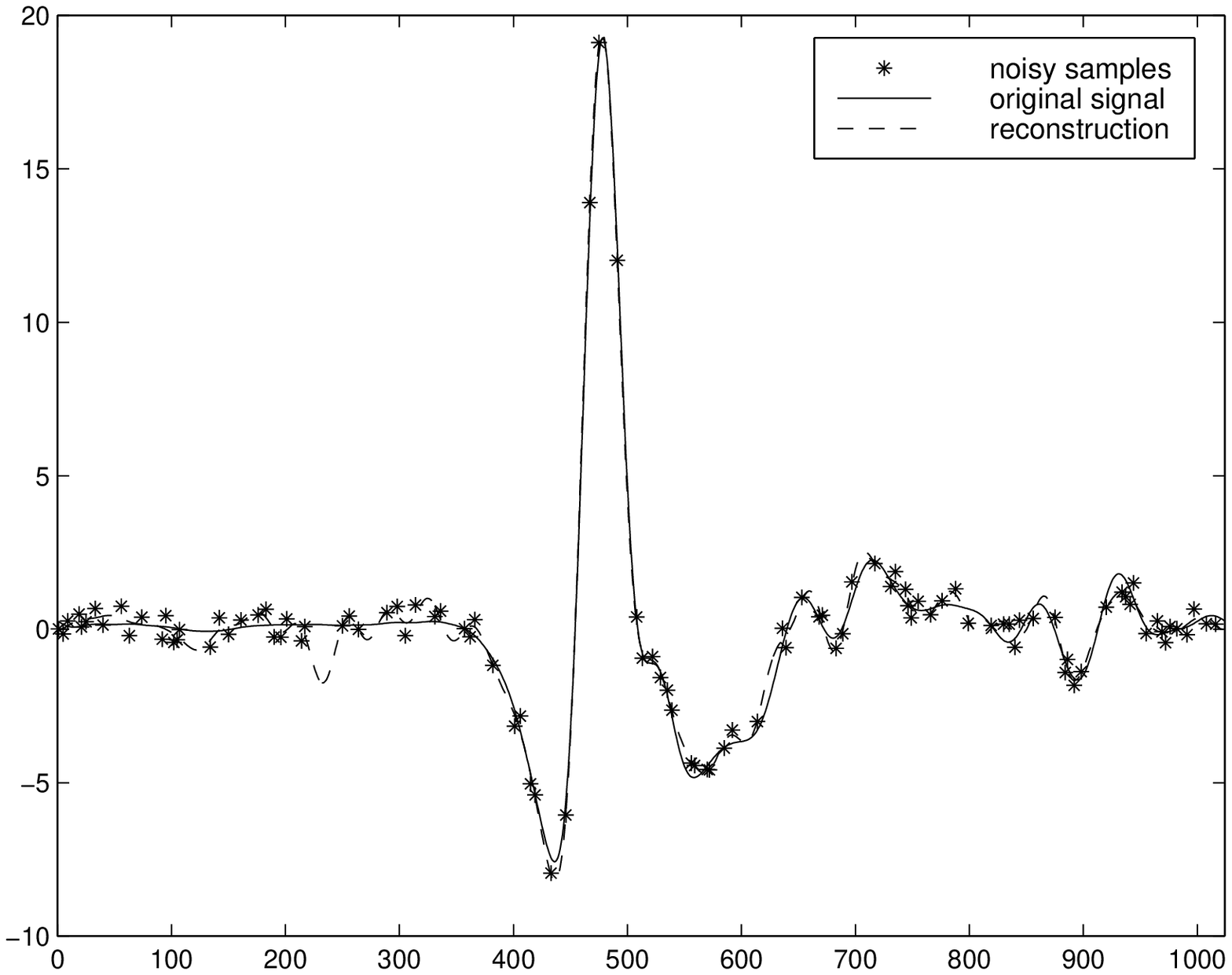,width=55mm,height=40mm}}
\subfigure[Using a too small bandwidth,
error=0.4645.]{
\epsfig{file=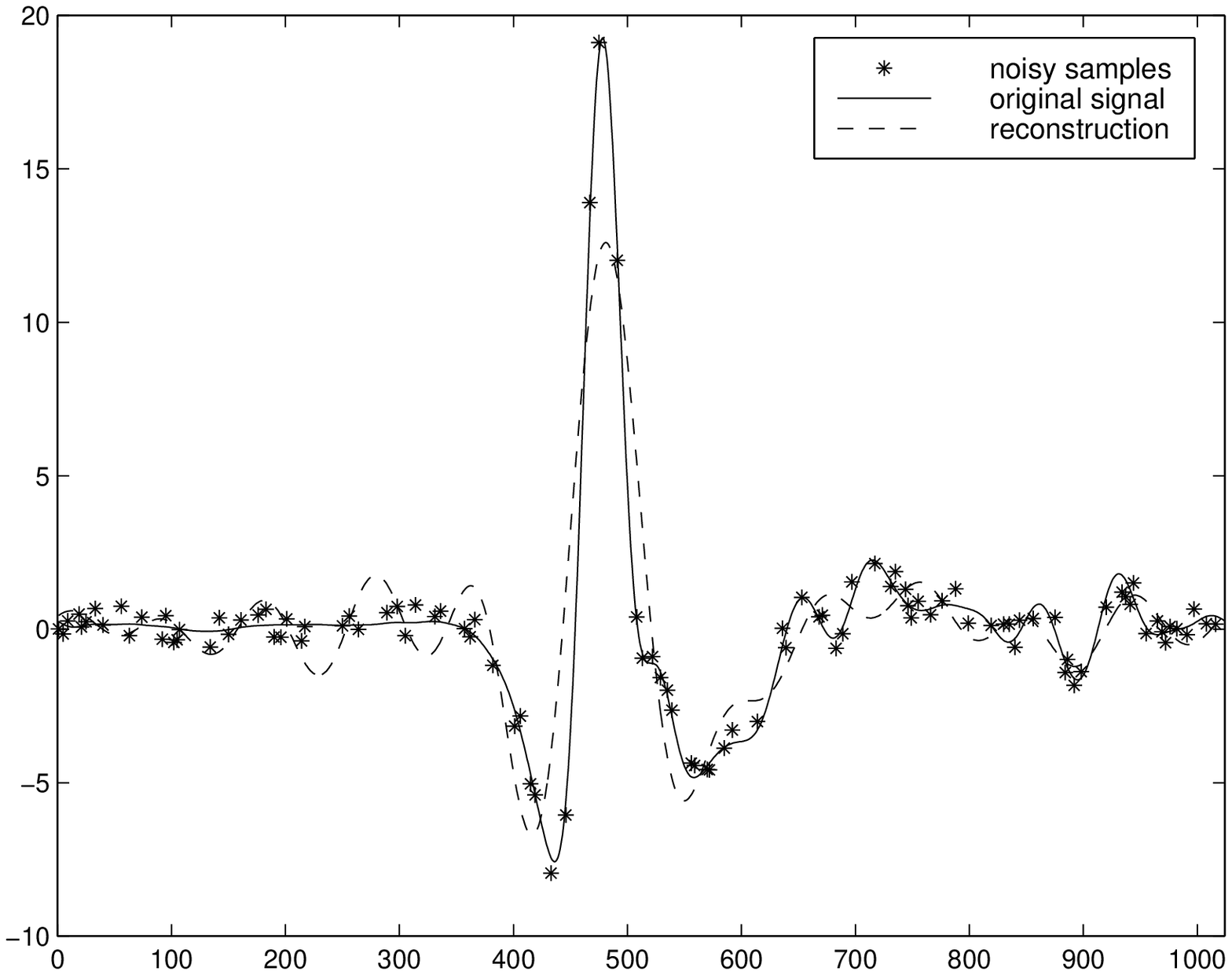,width=55mm,height=40mm}}\\
\subfigure[Using a too large bandwidth,
error = 0.2412.]{
\epsfig{file=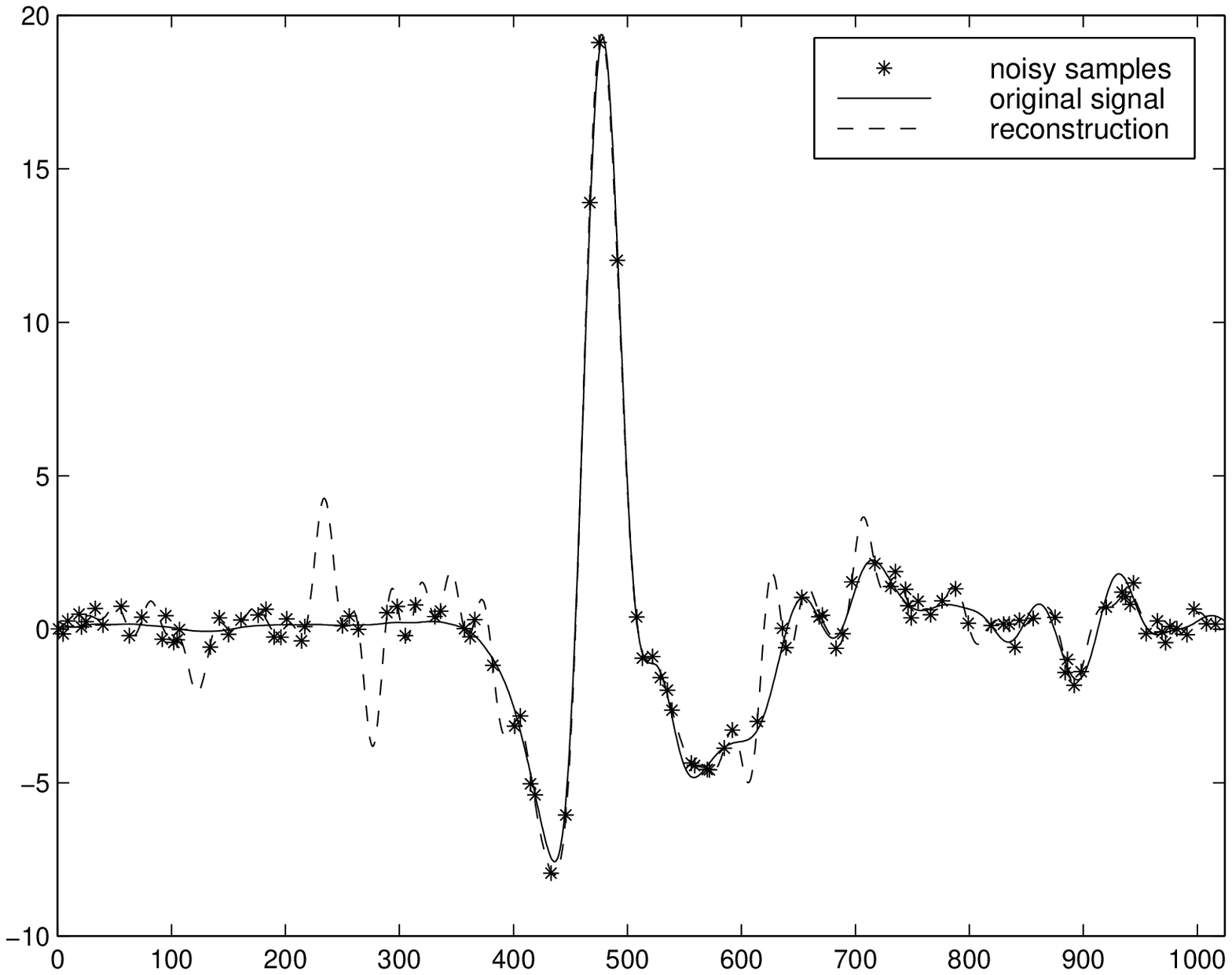,width=55mm,height=40mm}}
\subfigure[Multilevel algorithm, error=0.0959.]{
\epsfig{file=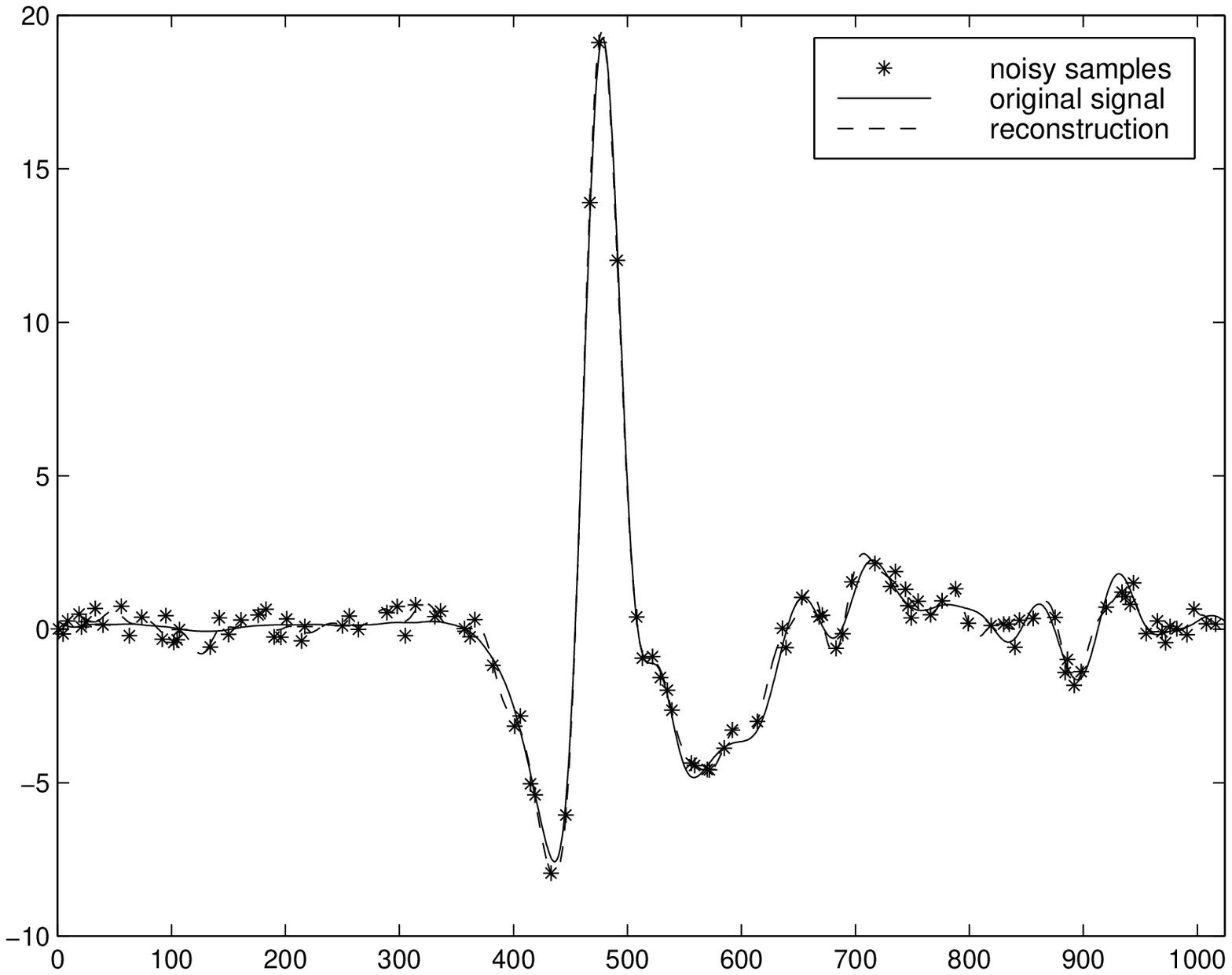,width=55mm,height=40mm}}
\caption{Example from spectroscopy -- comparison of reconstruction methods.}
\label{fig:spect}
\end{center}
\end{figure}

First we apply the truncated frame method with regularized SVD as
described in Section~\ref{ss:truncated}. We choose the truncation level
for the SVD via formula~\eqref{thresh}. This is the optimal
truncation level in this case, providing an approximation with least
squares error $0.0944$. Figure~\ref{fig:spect}(a) shows the reconstructed
signal together with the original signal and the noisy samples.
Without regularization we get a much worse ``reconstruction''
(which is not displayed).

We apply CG to the truncated frame method, as proposed in 
Section~\ref{ss:cgtrunc} with stopping criterion~\eqref{stopcg} (for
$\tau =1$). The algorithm terminates already after 3 iterations.
The reconstruction error is with $0.1097$ slightly higher than 
for truncated SVD (see also ~Figure~\ref{fig:spect}(b)), but the 
computational effort is much smaller.

Also Algorithm~\ref{th:act} (with $M=30$) terminates after 3 iterations.
The reconstruction is shown in~Figure~\ref{fig:spect}(c), the least squares 
error ($0.0876$) is slightly smaller than for the truncated frame
method, the computational effort is significantly smaller.

We also simulate the situation where the bandwidth is not known a priori
and demonstrate the importance of a good estimate of the bandwidth.
We apply Algorithm~\ref{th:act} using a too small
degree ($M = 11$) and a too high degree ($M = 40$). (We get
qualitatively the same results using the truncated frame method when
using a too small or too large bandwidth).
The approximations are shown in Figs.~\ref{fig:spect}(d) and (e),
The approximation errors are $0.4648$ and $0.2805$, respectively.
Now we apply the multilevel algorithm of Section~\ref{ss:ml}
which does not require any initial choice of the degree $M$.
The algorithm terminates at ``level'' $M=22$, the approximation
is displayed in Fig.~\ref{fig:spect}(f), the error is $0.0959$, thus
within the error bound $\delta$, as desired.
Hence without requiring explicit information about the bandwidth,
we are able to obtain the same accuracy as for the methods 
above.

\subsection{Approximation of geophysical potential fields}
\label{ss:geo}
Exploration geophysics relies on surveys of the Earth's magnetic field
for the detection of anomalies which reveal underlying geological features.
Geophysical potential field-data are generally observed at scattered
sampling points. Geoscientists, used to looking
at their measurements on maps or profiles and aiming at further
processing, therefore need a representation of the originally
irregularly spaced data at a regular grid. 

The reconstruction of a 2-D signal from its scattered data
is thus one of the first and crucial steps in geophysical
data analysis, and a number of practical constraints such as measurement
errors and the huge amount of data make the development of reliable
reconstruction methods a difficult task.

It is known that the Fourier transform of a geophysical potential field $f$ 
has decay $|\hat{f}(\omega)| = \ord (e^{-|\omega|})$. This rapid decay 
implies that $f$ can be very well approximated
by band-limited functions~\cite{RS98}. Since in general we may not 
know the (essential) bandwidth of $f$, we can use the
multilevel algorithm proposed in Section~\ref{ss:ml} to reconstruct $f$.

The multilevel algorithm also takes care of following problem.
Geophysical sampling sets are often highly anisotropic and large
gaps in the sampling geometry are very common. The large gaps
in the sampling set can make the reconstruction problem ill-conditioned
or even ill-posed. As outlined in Section~\ref{ss:regul} the multilevel 
algorithm iteratively determines the optimal bandwidth that 
balances the stability and accuracy of the solution.

\begin{figure}
\label{fig:geo}
\begin{center}
    \subfigure[Contour map of synthetic gravity anomaly, gravity is in
mGal.]{
    \epsfig{file=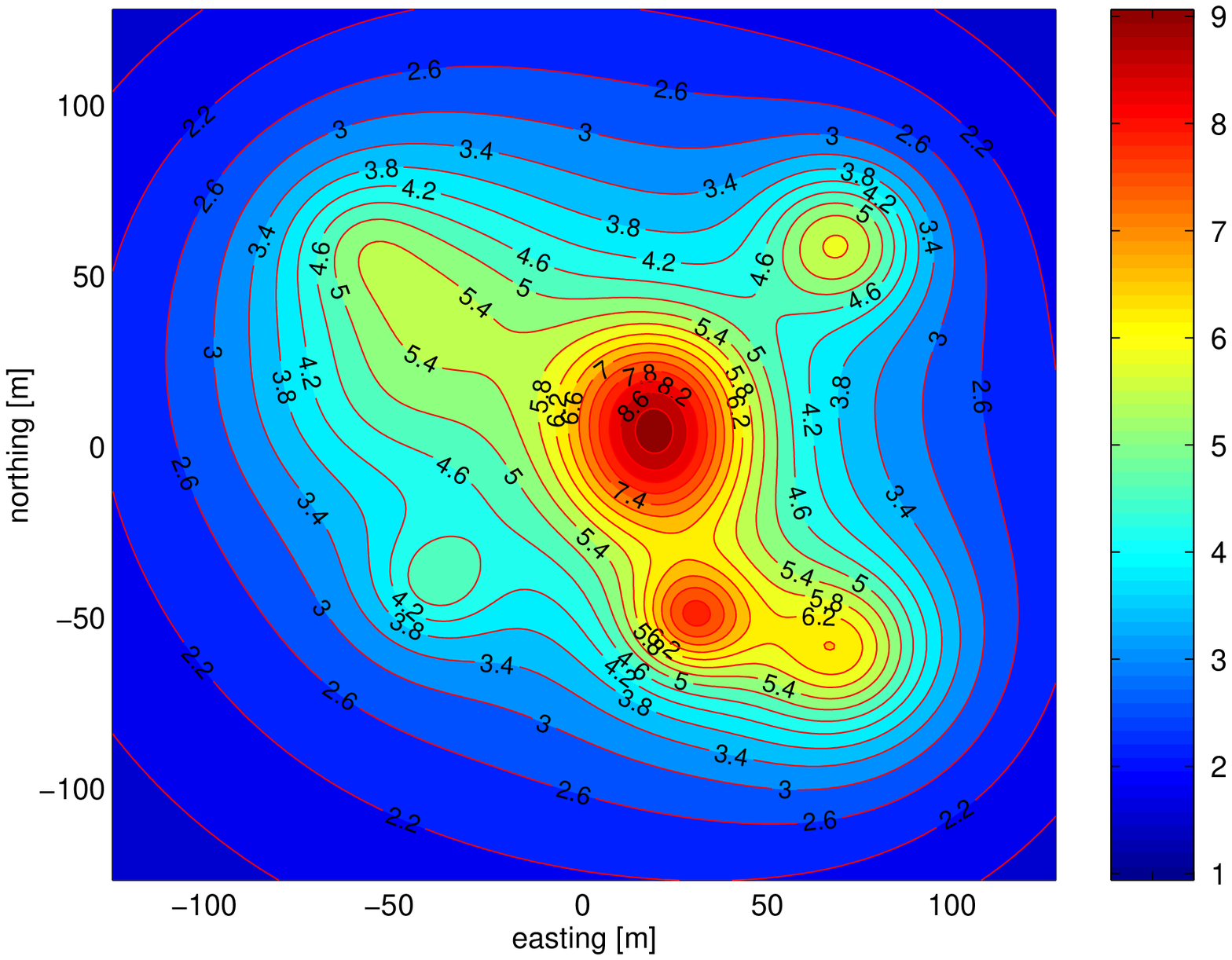,width=58mm,height=50mm}} \quad
    \subfigure[Sampling set and synthetic gravity anomaly.]{
    \epsfig{file=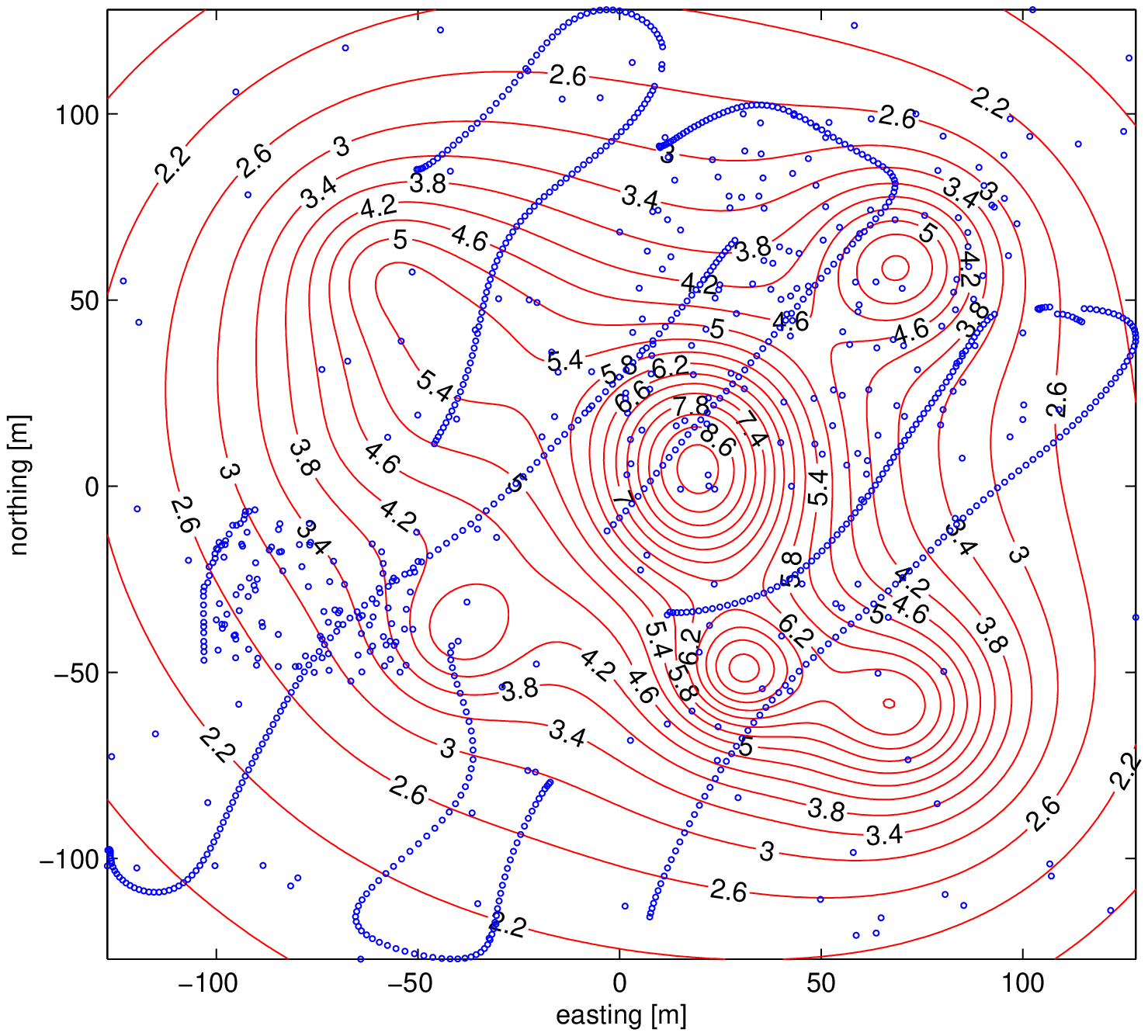,width=60mm,height=50mm}} \\
    \subfigure[Approximation by multi-level algorithm]{
    \epsfig{file=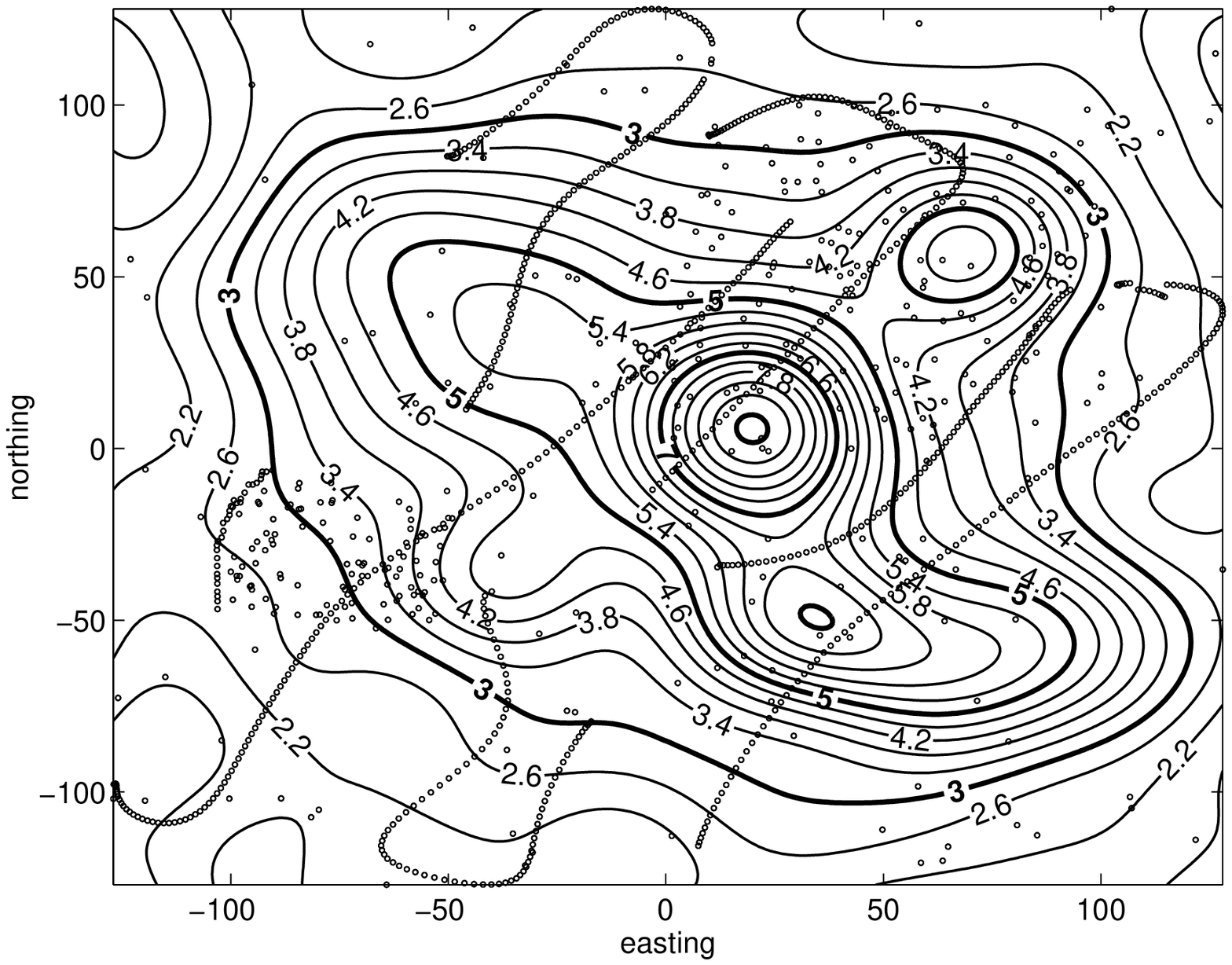,width=60mm}}\quad
    \subfigure[Error between approximation and actual anomaly.]{
    \epsfig{file=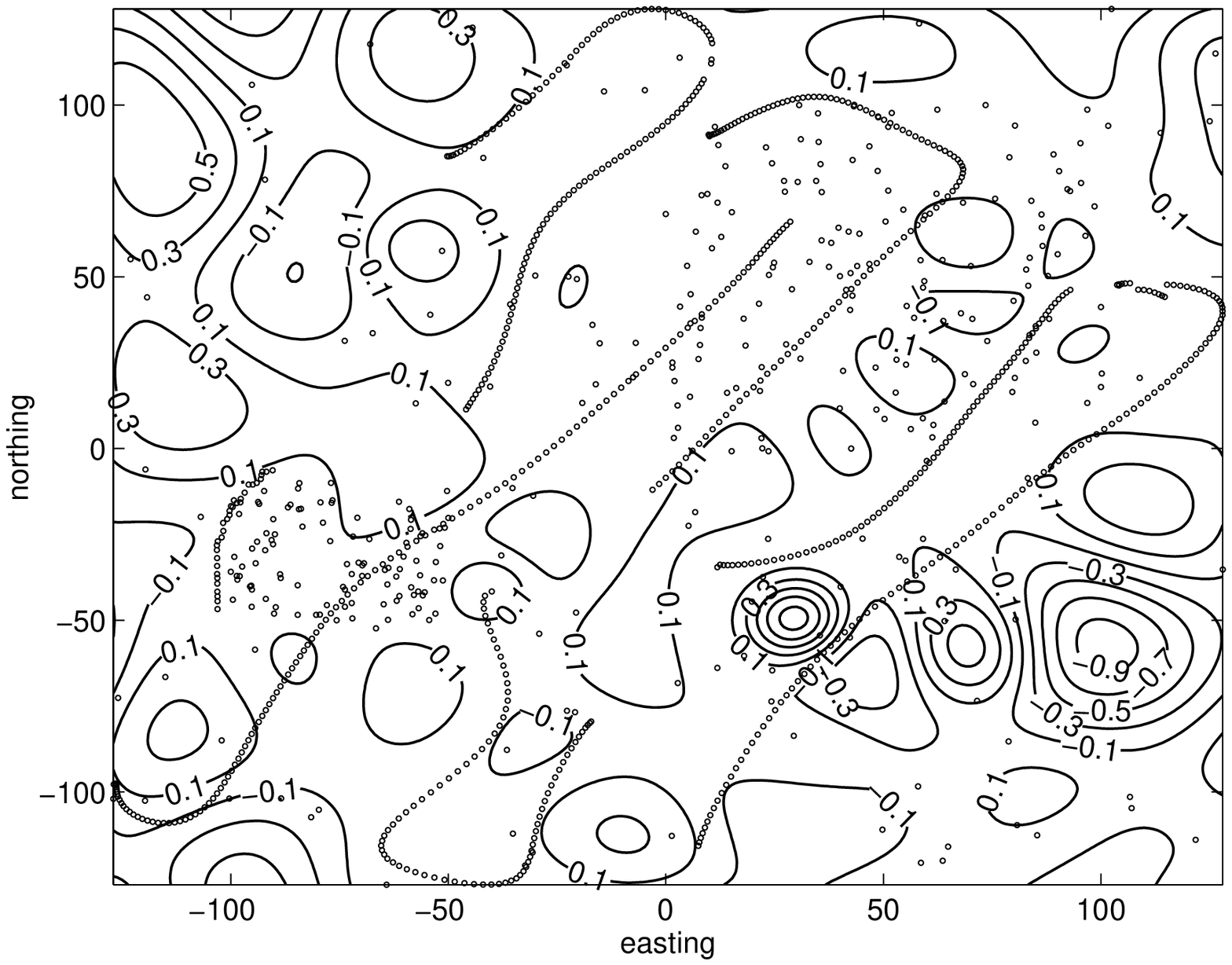,width=60mm}}
    \label{fig:geo2}
    \caption{Approximation of synthetic gravity anomaly from 1000
     non-uniformly spaced noisy samples by the multilevel algorithm of 
     Section~\ref{ss:ml}. The algorithm iteratively determines the
     optimal bandwidth (i.e.\ level) for the approximation.}
\end{center}
\end{figure}

Figure~\ref{fig:geo}(a) shows a synthetic gravitational anomaly $f$.
The spectrum of $f$ decays exponentially, thus
the anomaly can be well represented by a band-limited function,
using a ``cut-off-level'' of $|f(\omega)|  \le 0.01$ for the essential
bandwidth of $f$. 

We have sampled the signal at 1000 points $(u_j,v_j)$ and added 5\% random 
noise to the sampling values $f(u_j,v_j)$. The sampling geometry -- shown in 
Figure~\ref{fig:geo} as black dots -- exhibits several features one 
encounters frequently in exploration geophysics~\cite{RS98}. 
The essential bandwidth of $f$ would imply to choose a polynomial degree
of $M=12$ (i.e., $(2M+1)^2 = 625$ spectral coefficients).
With this choice of $M$ the corresponding block Toeplitz matrix $T_M$ 
would become ill-conditioned, making the reconstruction problem unstable.
As mentioned above, in practice we usually do not know the 
essential bandwidth of $f$. Hence we will not make use 
of this knowledge in order to approximate $f$.

We apply the multilevel method 
to reconstruct the signal, using only the sampling points
$\{(u_j,v_j)\}$, the samples $\{f^{\delta}(u_j,v_j)\}$ and the noise level
$\delta=0.05$ as a priori information. The algorithm
terminates at level $M=7$. The reconstruction is
displayed in Figure~\ref{fig:geo}(c), the error between the true
signal and the approximation is shown in Figure~\ref{fig:geo}(d).
The reconstruction error is $0.0517$ (or $0.193$ mGal), thus of the
same order as the data error, as desired.

\end{document}